\documentclass[11pt]{article}

\usepackage{fullpage}
\usepackage{authblk}
\usepackage{natbib}
\usepackage{algorithm}
\usepackage{algpseudocode}
\usepackage{amsmath,amsthm,amsfonts}
\usepackage{amsmath}
\usepackage[subnum]{cases}
\usepackage{tcolorbox}
\usepackage{mathrsfs}
\usepackage{amssymb}
\usepackage{color}
\usepackage{mathrsfs}
\usepackage{enumitem}
\usepackage{bm}
\usepackage{multirow}
\usepackage{booktabs}
\usepackage{makecell}
\usepackage{graphicx}
\usepackage{subcaption}
\usepackage{comment}
\usepackage{cases}
\usepackage{appendix}
\usepackage{tikz}
\usetikzlibrary{arrows,shapes}
\usepackage[colorlinks= true, linkcolor=red, citecolor=blue, urlcolor=black]{hyperref}
\usepackage{cleveref}
\usepackage{marginnote}
\usepackage{diagbox} 

\numberwithin{equation}{section}

\theoremstyle{definition}

\newtheorem{theorem}{Theorem}[section]

\newtheorem{corollary}[theorem]{Corollary}

\usepackage{booktabs}
\usepackage{tikz}
\usepackage{pgfplots}
\usetikzlibrary{calc,intersections}

\usetikzlibrary{shapes.geometric, arrows, positioning}

\tikzset{
	mybox/.style  = {draw, rectangle, minimum width=4cm, minimum height=0.8cm, text centered, text width=4.4cm,   
		font=\normalsize},
	box/.style  = {draw, rectangle, minimum width=2.0cm, minimum height=0.6cm, text centered, text width=3.0cm,   
		font=\normalsize},
	myarrow/.style = {line width=0.2pt, draw=black, -triangle 60, postaction={draw, line width=0.2pt, shorten >=10pt,-}}
}

\tikzstyle{arrow} = [->, >=stealth, -triangle 60]

\allowdisplaybreaks

\makeatletter
\newcommand{\leqnomode}{\tagsleft@true}
\newcommand{\reqnomode}{\tagsleft@false}
\makeatother

\begin{document}

\title{Gradient Norm Minimization of Nesterov Acceleration: $o(1/k^3)$}


\author[1,2]{Shuo Chen} 
\author[1,2]{Bin Shi}
\author[1,2]{Ya-xiang Yuan}
\affil[1]{State Key Laboratory of Scientific and Engineering Computing, Academy of Mathematics and Systems Science, Chinese Academy of Sciences, Beijing 100190, China \thanks{Email: \url{{chenshuo, shi bin, yyx}@lsec.cc.ac.cn}} }
\affil[2]{University of Chinese Academy of Sciences, Beijing 100049, China}

\date\today

\maketitle

\begin{abstract}
In the history of first-order algorithms, Nesterov's accelerated gradient descent (\texttt{NAG}) is one of the milestones. However, the cause of the acceleration has been a mystery for a long time. It has not been revealed with the existence of gradient correction until the high-resolution differential equation framework proposed in~\citep{shi2021understanding}. In this paper, we continue to investigate the acceleration phenomenon. First, we provide a significantly simplified proof based on precise observation and a tighter inequality for $L$-smooth functions. Then, a new implicit-velocity high-resolution differential equation framework, as well as the corresponding implicit-velocity version of phase-space representation and Lyapunov function, is proposed to investigate the convergence behavior of the iterative sequence $\{x_k\}_{k=0}^{\infty}$ of~\texttt{NAG}. Furthermore, from two kinds of phase-space representations, we find that the role played by gradient correction is equivalent to that by velocity included implicitly in the gradient, where the only difference comes from the iterative sequence $\{y_{k}\}_{k=0}^{\infty}$ replaced by $\{x_k\}_{k=0}^{\infty}$. Finally, for the open question of whether the gradient norm minimization  of~\texttt{NAG} has a faster rate $o(1/k^3)$, we figure out a positive answer with its proof. Meanwhile, a faster rate of objective value minimization $o(1/k^2)$ is shown for the case $r > 2$.
\end{abstract}

%
%
%
%

\section{Introduction}
\label{sec: intro}

Gradient-based optimization, as the workhorse algorithm, has been powering recent developments in statistical machine learning since typically the first-order information can be computed cheaply and available in practice.  Recall the unconstrained minimization problems with $f$ being a $L$-smooth convex function, 
\[
\min_{x \in \mathbb{R}^d} f(x),
\]
with $f$ being a $L$-smooth convex function, a milestone in the development works of gradient-based algorithms is Nesterov's accelerated gradient descent (\texttt{NAG}),
\[
\label{eqn: nesterov}
\left\{\begin{aligned}
& x_k = y_{k-1} - s \nabla f(y_{k-1}) \\
& y_{k} = x_{k} + \frac{k-1}{k+r}(x_{k} - x_{k-1}),
\end{aligned} \right.
\]
with $x_0 = y_0 \in \mathbb{R}^d$, of which the case $r=2$ is originally proposed by~\citet{nesterov1983method}. Then,~\citet{nesterov2003introductory} also proposes an accelerated gradient scheme for  $\mu$-strongly convex objective functions with a slight change in the momentum coefficients. For any step size $0 < s \leq 1/L$, if the objective function is $\mu$-strongly convex, \texttt{NAG} converges with a linear rate as  
\[
f(x_k) - f(x^\star) \leq O\left( (1 - \sqrt{\mu s})^k \right),
\]
which essentially improves the rate $ O\left( (1 - \mu s)^k \right)$ achieved by the standard gradient descent; meanwhile, if the objective function is convex,  \texttt{NAG} achieves an inverse quadratic rate as  
\[
f(x_k) - f(x^\star) \leq O\left( \frac{1}{sk^2} \right),
\]
while the convergence rate generated by the standard gradient descent is only $O(1/sk)$. 

Although the accelerated convergence rates have been obtained by~\citet{nesterov1983method,nesterov2003introductory}, the cause leading to acceleration cannot be found with the technique --- ``Estimate Sequence'' originally developed by Nesterov himself. Recently,~\citet{shi2021understanding} expands~\texttt{NAG} with the single sequence $\{y_{k}\}_{k=0}^{\infty}$ as
\begin{equation}
\label{eqn: nag-sy}
y_{k+1} = y_k + \frac{k}{k+r+1}\left( y_{k} - y_{k-1} \right) - s \nabla f(y_k) - \frac{k}{k+r+1} \cdot s \underbrace{\left( \nabla f(y_{k}) - \nabla f(y_{k-1}) \right)}_{\textbf{gradient correction}},
\end{equation}
where we find the existence of gradient correction, which is different from Polyak's heavy-ball method (or momentum method), so we call~\eqref{eqn: nag-sy} \textbf{gradient-correction scheme of \texttt{NAG}}. The high-resolution differential equation framework with the corresponding phase-space representation and Lyapunov function techniques is proposed in~\citep{shi2021understanding}, which confirms the cause leading to acceleration is exactly the existence of gradient correction. However, if we represent~\texttt{NAG} with the following single sequence $\{x_{k}\}_{k=0}^{\infty}$ as  
\begin{equation}
\label{eqn: nag-sx}
x_{k+1} = x_{k} + \frac{k-1}{k+r} \left(x_{k} - x_{k-1} \right) - s \nabla f\bigg(x_{k} + \frac{k-1}{k+r} \underbrace{\left(x_{k} - x_{k-1} \right)}_{\textbf{implicit velocity}} \bigg),
\end{equation}
where the gradient correction term in the scheme~\eqref{eqn: nag-sy} is replaced by the position correction term included implicitly in the gradient in the scheme~\eqref{eqn: nag-sx}. In the following sections, we will know that position correction is defined as velocity,  so we call~\eqref{eqn: nag-sx} \textbf{implicit-velocity scheme of~\texttt{NAG}}. The viewpoint that the velocity included implicitly in the gradient will generate the acceleration in the continuous differential equation has been noticed in~\citep{muehlebach2019dynamical, alecsa2021extension, adly2022convergence}. However, a large gap exists in obtaining the convergence rate for the discrete algorithms, at least for the implicit-velocity scheme of~\texttt{NAG}~\eqref{eqn: nag-sx} itself. Hence, we mention it as the following question. 
\begin{tcolorbox}
\begin{center}
\textbf{Question:}  Can the same convergence rate be obtained from the implicit-velocity scheme of~\texttt{NAG}~\eqref{eqn: nag-sx}?
\end{center}
\end{tcolorbox}


  In this paper, based on the high-resolution ODE framework as well as the new techniques developed in~\citep{shi2021understanding} --- phase-space representation and Lyapunov function, we will investigate further the general case of~\texttt{NAG} ($r \geq 2$) from both the gradient-correction scheme~\eqref{eqn: nag-sy} and the implicit-velocity scheme~\eqref{eqn: nag-sx}.

\subsection{Overview of contributions}
\label{subsec: overview-contributions}

In this paper, we study further the convergence rate of~\texttt{NAG} based on the high-resolution differential equation framework proposed in~\citep{shi2021understanding} and make the following contributions.

\paragraph{Simplification of the proof in~\citep{shi2021understanding}}  
We revisit the proof about the convergence rate of~\texttt{NAG} in~\citep{shi2021understanding}. Based on precise observation, the forward calculation of the difference between the mixed terms of the Lyapunov function can lead to some complex terms being canceled out. With the tighter inequality for the $L$-smooth function, the proof together with the coefficients in the convergence rate for~\texttt{NAG} is highly simplified. Moreover, we extend the proof to the general case $r \geq 2$. 

\paragraph{Proof from the implicit-velocity scheme of~\texttt{NAG} } Along the way of high-resolution differential equations~\citep{shi2021understanding}, we start to consider the~\texttt{NAG} from the implicit scheme~\eqref{eqn: nag-sx}. The corresponding implicit-velocity high-resolution differential equation with its continuous convergence rate is derived. The phase-space representation for the implicit-velocity scheme~\eqref{eqn: nag-sx} is implicit-explicit; on the contrary, that is explicit-implicit for the gradient-correction scheme~\eqref{eqn: nag-sy}. Although the Lyapunov function for the implicit-velocity scheme~\eqref{eqn: nag-sx} is constructed along the way different from the gradient-correction scheme~\eqref{eqn: nag-sy}, the two Lyapunov functions are essentially identical. Finally, based on the convergence rates obtained, we find the Lyapunov function based on the high-resolution resolution differential equation framework with its phase-space representation is tailor-made for~\texttt{NAG}.


\paragraph{New acceleration phenomena}
For the convergence rates of~\texttt{NAG}, we first answer the open problem proposed in~\citep{shi2021understanding} as
\[
\min_{0\leq i \leq k}\| \nabla f(y_i) \|^2 \leq o\left( \frac{1}{k^3} \right),
\]
which is based on precise observation and calculation of the sum of infinite series. Meanwhile, we also point out similarly that the convergence rate for the objective proposed in~\citep{su2016differential}  can be improved to
\[
\min_{0\leq i \leq k}f(y_i) - f(x^\star) \leq o\left( \frac{1}{k^2} \right),
\]
which is also available in practice and far more straightforward than the derivation in~\citep{attouch2016rate}.

\subsection{Organization and notations}
\label{subsec: notation-organization}
The remainder of the paper is organized as follows. In~\Cref{sec: related-works}, some related research works are described. The simplified proof of general~\texttt{NAG} ($r \geq 2$) is provided in~\Cref{sec: simplification}. We investigate the convergence rate of general~\texttt{NAG} from the iterative sequence $\{x_k\}_{k=0}^{\infty}$ based on the implicit-velocity version of  high-resolution differential equation framework in~\Cref{sec: implicit-proof}. In~\Cref{sec: discussion}, we discuss the similarities and differences between the two Lyapunov functions and propose further research directions. 

In this paper, we follow the notation of~\citep{shi2021understanding}. Let $\mathcal{F}_L^1(\mathbb{R}^d)$ be the class of $L$-smooth convex functions defined on $\mathbb{R}^d$; that is, $f \in \mathcal{F}_L^1$ if $f(y) \geq f(x) + \langle \nabla f(x), y - x\rangle$ for all $x, y \in \mathbb{R}^d$ and its gradient is $L$-Lipschitz continuous in the sense that $\| \nabla f(x) - \nabla f(y) \| \leq L\|x - y\|$, where $\|\cdot\|$ denotes the standard Euclidean norm and $L > 0$ is the Lipschitz constant. And $x^\star$ denotes the minimizer of the convex objective $f$. 

\section{Related works}
\label{sec: related-works}

The history of modern accelerated gradient descent methods dates from~\citep{gelfand1962some}, which first proposes a two-step Ravine method only with the gradient information.~\citet{polyak1964some} provides his heavy-ball method (or momentum method) based on the invariant manifold theorem from the field of dynamical systems, which locally accelerates the convergence rate compared with the standard gradient descent. The milestone works are due to~\citet{nesterov1983method,nesterov2003introductory} for his accelerated gradient methods, which improve the convergence rate globally for both the $\mu$-strongly convex and convex functions. Because there are so many complex algebraic tricks in Nesterov's own ``estimate sequence'' techniques, some open questions are proposed to find the cause that leads to the acceleration~\citep{bubeck2015convex}.  Moreover, a significant generalization with broad applications to nonsmooth objective functions is provided in~\citep{beck2009fast}.

Recently,~\citet{su2016differential} proposed a low-resolution differential equation framework to model the~\texttt{NAG} and used the Lyapunov function's techniques to investigate the convergence rate. A series of works based on the Lagrangian and Hamiltonian frameworks are used to
generate a large class of continuous-time differential equations for a unified treatment of accelerated gradient-based methods~\citep{wibisono2016variational,wilson2021lyapunov}. Indeed, from a single ``Bregman Lagrangian'', the~\texttt{NAG} is extended to non-Euclidean settings, such as mirror descent and accelerated higher-order gradient methods~\citep{wibisono2016variational}. An equivalence between the estimate sequence technique and the Lyapunov function technique is established in~\citep{wilson2021lyapunov}, which further strengthens the connection between differential equations and discrete algorithms. A detailed review is shown in~\citep{jordan2018dynamical}.

Another venerable line of works combining the low-resolution differential equation and the continuous limit of the Newton method to design and analyze new algorithms is proposed in~\citep{alvarez2000minimizing, attouch2012second, attouch2014dynamical, attouch2016fast, attouch2016rate}, of which the terminology is called inertial dynamics with a Hessian-driven term. Although the inertial dynamics with a Hessian-driven term resembles closely with the high-resolution differential equations in~\citep{shi2021understanding}, it is important to note that the Hessian-driven terms are from the second-order information of Newton's method~\citep{attouch2014dynamical}, while the gradient correction entirely relies on the first-order information of~\texttt{NAG}. Also, there is an exciting work that achieves the gradient norm acceleration $O(1/k^3)$ by introducing an additional sequence of iterates and a more aggressive step-size policy to the original~\texttt{NAG}~\citep{ghadimi2016accelerated}.

For the gradient-based algorithms, the stochastic setting is more valuable and practical. The underdamped Langevin diffusion with the log-concave target distribution is investigated, and an MCMC algorithm based on its discretization is provided in~\cite{cheng2018underdamped}. A more comprehensive work about overdamped and underdamped Langevin MCMC with the theoretical upper bounds on the number of steps is shown in~\citep{cheng2018sharp}. Furthermore, based on the Langevin dynamics with its macroscopic evolving differential equation as well as the spectral structure of Schr\"odinger operators,~\citet{shi2020learning} proposes a profound theoretical analysis for the stochastic gradient descent on the nonconvex objectives and points out the effect of the learning rate. Afterward, the theoretical analysis for the hyperparameters in stochastic gradient descent with momentum is demonstrated in~\citep{shi2021hyperparameters}, which is based on the more complicated kinetic Fokker-Planck equation.


\section{A simplified proof and generalization }
\label{sec: simplification}

In this section, starting from the gradient-correction scheme~\eqref{eqn: nag-sy} in~\citep{shi2021understanding}, we revisit the original proof about the convergence rate of~\texttt{NAG} with $r=2$ . Then, based on the observation of a more tight inequality, we propose a simplified proof and generalize it to the cases $r \geq 2$.

First, with the velocity iterates $v_{k-1} = (y_{k} - y_{k-1})/\sqrt{s}$, we look back the phase-space representation of implicit-velocity scheme of~\texttt{NAG} as 
\begin{equation}
\label{eqn:phase-nag}
\left\{ \begin{aligned}
        & y_{k} - y_{k-1} = \sqrt{s}v_{k-1} \\
        & v_{k} - v_{k-1} = - \frac{r+1}{k}v_{k} -  \left( 1 + \frac{r+1}{k}\right)\sqrt{s}\nabla f\left( y_{k}  \right)  - \sqrt{s}\left(\nabla f(y_{k}) - \nabla f(y_{k-1})\right),
        \end{aligned}\right.  
\end{equation}
with any initial $y_0$ and $v_0 = -\sqrt{s} \nabla f(y_0)$. Here we can find in the phase-space representation~\eqref{eqn:phase-nag}, the first line for the position sequence $\{y_k\}_{k=0}^{\infty}$ is an explicit scheme, while the second line for the velocity sequence $\{v_{k}\}_{k=0}^{\infty}$ is an implicit scheme. Moreover, the gradient-correction term, $\sqrt{s}\left( \nabla f(y_k) - \nabla f(y_{k-1})  \right)$, is shown in the velocity iteration, so  we call~\eqref{eqn:phase-nag}~\textbf{gradient-correction phase-space representation of~\texttt{NAG}}. Recall the Lyapunov function for \texttt{NAG} with $r=2$ based on the framework of high-resolution differential equations is constructed in~~\citep[(4.41)]{shi2021understanding}, 
\begin{equation}
\label{eqn:lypunov-441}
\mathcal{E}(k) = s(k+1)(k+3)\left( f(y_{k}) - f(x^\star) \right) + \frac12 \| \sqrt{s}(k + 1)v_{k} + 2(y_{k+1} - x^\star) + s(k + 1)\nabla f(y_{k}) \|^2,
\end{equation}
and then it is deduced in~\citep[Theorem 6]{shi2021understanding} that the iterative sequence $\{y_k\}_{k=0}^{\infty}$ converges as 
\begin{equation}
\label{eqn:rate-shi2021}
 f(y_{k}) - f(x^\star) \leq \frac{119\|y_0 - x^\star\|^2}{s(k+1)^2},\qquad  \min_{0\leq i\leq k}\| \nabla f(y_i) \|^2 \leq \frac{8568\|y_0 - x^\star\|^{2}}{s^2(k+1)^{3}},
\end{equation}
for any step size $s \leq 1/(3L)$.  Although the second inequality of~\eqref{eqn:rate-shi2021} shows the acceleration of gradient norm minimization, a simple glance tells us that the coefficients of two inequalities~\eqref{eqn:rate-shi2021} are too large. Compared with~\citep[Theorem 2.2.2]{nesterov2003introductory}, the step size is required to shrink by a factor of $3$, and the coefficient of the first inequality is also overlarge. Moreover, to estimate the first several iterates,~\citet[Section 4.2 and Section C]{shi2021understanding} uses a plenary of gradient Lipschitz inequality, which leads to the coefficient being overlarge.  

For the $L$-smooth function, the estimate of difference between $f(y_{k+1})$ and $f(x^\star)$ is used in ~\citep[Section C.4.4]{shi2021understanding} as
\[
f(y_{k+1}) - f(x^\star)  \leq \left\langle \nabla f(y_{k+1}), y_{k+1} -x^\star \right\rangle. 
\]
As a matter of simple fact, the estimate can be bounded by a tighter inequality as 
\begin{equation}
\label{eqn:lipschitz-grad}
f(y_{k+1}) - f(x^\star)  \leq \left\langle \nabla f(y_{k+1}), y_{k+1} -x^\star \right\rangle - \frac{1}{2L} \| \nabla f(y_{k+1}) \|^2.
\end{equation}
Before we proceed to the proof, we move one space backward for the general Lyapunov function~\eqref{eqn:lypunov-441} as 
\begin{equation}
\label{eqn: lypunov-positive}
\mathcal{E}(k) = sk(k+r)\left( f(y_{k-1}) - f(x^\star) \right) + \frac12 \| \sqrt{s}kv_{k-1} + r(y_k - x^\star) + sk\nabla f(y_{k-1}) \|^2.
\end{equation}
With the tighter inequality~\eqref{eqn:lipschitz-grad}, we conclude the result with the following theorem and simple proof.

\begin{theorem}
\label{thm:nag}
Let $f \in \mathcal{F}_{L}^{1}(\mathbb{R}^d)$. For any step size $0 < s \leq 1/L$, the iterates $\{y_{k}\}_{k=0}^{\infty}$ generated by \texttt{NAG} obey
\begin{equation}
\label{eqn:y-rate-nag}
 f(y_{k}) - f(x^\star) \leq \frac{r^2\|x_0 - x^\star\|^2}{2s(k+1)(k+r+1)},\quad  \min_{0\leq i \leq k} \big\| \nabla f(y_{i}) \big\|^2 \leq \frac{6r^2 \|x_0 - x^\star\|^2}{s^2(k+1)(k+2)(2k+3r+3)}; 
\end{equation}
furthermore, the iterative sequence $\{y_{k}\}_{k=0}^{\infty}$ also satisfy
\begin{equation}
\label{eqn:y-rate-nag-o}
\lim_{k \rightarrow \infty} \left( k^3 \min_{0\leq i \leq k} \big\| \nabla f(y_{i}) \big\|^2 \right) = 0.
\end{equation}
If $r > 2$, then the iterative sequence $\{y_{k}\}_{k=0}^{\infty}$ also satisfy
\begin{equation}
\label{eqn:y-rate-nagmore3}
\lim_{k \rightarrow \infty} \left[ k^2 \left( \min_{0\leq i \leq k}  f(y_{i}) - f(x^\star) \right)  \right] = 0.
\end{equation}
%

\end{theorem}

According to~\Cref{thm:nag}, the coefficients of the convergence rates in~\eqref{eqn:y-rate-nag} is simplified greatly when $r=2$, compared with the convergence rates~\eqref{eqn:rate-shi2021} deduced in~\citep[Theorem 6]{shi2021understanding}. 

\begin{proof}[Proof of~\Cref{thm:nag}]
Before starting the proof, we first reformulate the second iterative equality of the explicit-implicit gradient-correction phase-space representation~\eqref{eqn:phase-nag} equivalently as
\begin{equation}
\label{eqn:phase-nag2}
(k+r+1)v_{k} - kv_{k-1} + \sqrt{s}\left[(k+r+1)\nabla f(y_k) - k\nabla f(y_{k-1}) \right] =  - k\sqrt{s}\nabla f(y_k).
\end{equation}
For convenience, to make the proof clear, we separate it into three steps.
\begin{itemize}
\item[\textbf{(1)}] Before we proceed to the next for obtaining further the difference between the Lyapunov functions~\eqref{eqn: lypunov-positive},  $\mathcal{E}(k+1)$ and $\mathcal{E}(k)$,  the key point here is to compute the difference between the second terms, $\sqrt{s}(k+1)v_{k} + r(y_{k+1} - y^\star) + s(k+1)\nabla f(y_{k})$ and $\sqrt{s}kv_{k-1} + r(y_{k} - y^\star) + sk\nabla f(y_{k-1})$, as 
\begin{align*}
      & \left[\sqrt{s}(k+1)v_{k} + r(y_{k+1} - x^\star) + s(k+1)\nabla f(y_{k}) \right] - \left[\sqrt{s}kv_{k-1} + r(y_{k} - x^\star) + sk\nabla f(y_{k-1}) \right] \\
  =   & \sqrt{s}(k+r+1)v_{k} - \sqrt{s}kv_{k-1} + s(k+1)\nabla f(y_{k}) - sk\nabla f(y_{k-1}). 
\end{align*}
Then, with the equivalent form of the explicit-implicit gradient-correction phase-space representation~\eqref{eqn:phase-nag2}, the difference between the second terms can be calculated as
\begin{align}
      & \left[\sqrt{s}(k+1)v_{k} + r(y_{k+1} - x^\star) + s(k+1)\nabla f(y_{k}) \right] - \left[\sqrt{s}kv_{k-1} + r(y_{k} - x^\star) + sk\nabla f(y_{k-1}) \right] \nonumber \\
  =   & - s(k+r)\nabla f(y_{k}). \label{eqn:difference-a}
\end{align}

\item[\textbf{(2)}] Then, with the subscripts labeled, we calculate the difference between $\mathcal{E}(k+1)$ and $\mathcal{E}(k)$ as
\begin{align*}
\mathcal{E}&(k+1) - \mathcal{E}(k) \\
& =  \underbrace{ sk(k+r)\left( f(y_{k}) - f(y_{k-1}) \right) + s(2k+r+1)\left( f(y_{k}) - f(x^\star) \right) }_{\mathbf{I}}  + \frac{s^2(k+r)^2}{2}\big\| \nabla f(y_{k}) \big\|^2\\
                                  & \mathrel{\phantom{=}} - \underbrace{s(k+r)\left\langle \nabla f(y_{k}), \sqrt{s}kv_{k-1} + r(y_{k} - x^\star) + sk\nabla f(y_{k-1})\right\rangle}_{\mathbf{II}} \\
& =   \underbrace{sk(k+r)\left( f(y_{k}) - f(y_{k-1}) \right) }_{\mathbf{I}_1}+ \underbrace{s(2k+r+1)\left( f(y_{k}) - f(x^\star) \right)}_{\mathbf{I}_2}  + \frac{s^2(k+r)^2}{2}\big\| \nabla f(y_{k}) \big\|^2\\
                                  & \mathrel{\phantom{=}} - \underbrace{sk(k+r)\left\langle \nabla f(y_{k}), y_{k} - y_{k-1} \right\rangle }_{\mathbf{II}_1}- \underbrace{sr(k+r)\left\langle \nabla f(y_{k}),  y_{k} - x^\star\right\rangle}_{\mathbf{II}_2} - \underbrace{s^{2}k(k+r)\langle \nabla f(y_{k}), \nabla f(y_{k-1}) \rangle}_{\mathbf{II}_3} 
	\end{align*}

\item[\textbf{(3)}] Actually, we here can find the direct forward calculation will lead to a simplified proof. For any $f \in \mathcal{F}_{L}^{1}$, we can obtain the estimate as 
\begin{align*}
\mathbf{I}_1 - \mathbf{II}_1 & =        sk(k+r)\left( f(y_{k}) - f(y_{k-1}) \right) - sk(k+r)\left\langle \nabla f(y_{k}), y_{k} -  y_{k-1} \right\rangle \\
                             & \leq   - \frac{sk(k+r)}{2L} \| \nabla f(y_{k}) - \nabla f(y_{k-1}) \|^2; 
\end{align*}
moreover, with the tighter inequality~\eqref{eqn:lipschitz-grad}, the difference between $\mathbf{I}_2$ and $\mathbf{II}_2$ can be computed as
\begin{align*}
\mathbf{I}_2 - \mathbf{II}_2 & =   s(2k+r+1)\left( f(y_{k}) - f(x^\star) \right) - sr(k+r)\left\langle \nabla f(y_{k}),  y_{k} - x^\star\right\rangle \\
                             & \leq  - \frac{sr(k+r)}{2L} \| \nabla f(y_{k}) \|^2 - s \left[ (r-2)k + r^2 - r - 1  \right]\left( f(y_{k}) - f(x^\star) \right).
\end{align*}
Hence, the difference between $\mathcal{E}(k+1)$ and $\mathcal{E}(k)$ can be obtained as
\begin{align*}
\mathcal{E}(k+1) - \mathcal{E}(k) & =  - \frac{sk(k+r)}{2L} \| \nabla f(y_{k}) - \nabla f(y_{k-1}) \|^2  - \frac{sr(k+r)}{2L} \| \nabla f(y_{k}) \|^2 \\
                                  & \mathrel{\phantom{=}} + \frac{s^2(k+r)^2}{2}\big\| \nabla f(y_{k}) \big\|^2  - s^{2}k(k+r)\langle \nabla f(y_{k}), \nabla f(y_{k-1}) \rangle \\
                                  & \mathrel{\phantom{=}} - s \left[ (r-2)k + r^2 - r - 1  \right]\left( f(y_{k}) - f(x^\star) \right).
\end{align*}
When $0 < s \leq 1/L $, the difference between $\mathcal{E}(k+1)$ and $\mathcal{E}(k)$ satisfies
\begin{align}
	\label{eqn:lyapunov-monotone}
	\nonumber\mathcal{E}(k+1) - \mathcal{E}(k) \leq & - \frac{s^2k(k+r)}{2} \big\| \nabla f(y_{k-1}) \big\|^2 \\
	& - s \left[ (r-2)k + r^2 - r - 1  \right]\left( f(y_{k}) - f(x^\star) \right).
\end{align}
Hence, we can obtain the convergence rates~\eqref{eqn:y-rate-nag} with the derivation facts above that the Lyapunov function $\mathcal{E}(k)$ is decreasing monotonically, and the following inequality is satisfied as
\[
\frac12\sum_{i=0}^{\infty} s^2 (i+1)(i+r+1)\big\| \nabla f(y_{i}) \big\|^2 + \sum_{i=0}^{\infty} s \left[ (r-2)i + r^2 - r - 1  \right]\left( f(y_{i}) - f(x^\star) \right) \leq \frac{r^2 \|x_0 - x^\star\|^2}{2}.
\]
The limitation~\eqref{eqn:y-rate-nag-o} therefore is derived as
\[
0 \leq \frac{7}{24} \lim_{k \rightarrow \infty} \left( k^3 \min_{0\leq i \leq k} \big\| \nabla f(y_{i}) \big\|^2 \right) \leq \lim_{k \rightarrow \infty}\sum_{i= \left \lfloor \frac{k}{2}\right \rfloor}^{k} (i+1)(i+r+1)\big\| \nabla f(y_{i}) \big\|^2 = 0.
\]
Furthermore, when $r > 2$, we complete the derviation of the limitation~\eqref{eqn:y-rate-nagmore3} as
\[
0  \leq \frac{3}{8} \lim_{k \rightarrow \infty} \left[ k^2 \left( \min_{0\leq i \leq k}  f(y_{i}) - f(x^\star) \right)  \right]
  \leq \sum_{i=\left \lfloor \frac{k}{2}\right \rfloor}^{k}  \left[ i + \frac{r^2 - r - 1}{r-2}  \right]\left( f(y_{i}) - f(x^\star) \right) = 0.
\]

\end{itemize}
\end{proof}

For any $f\in \mathcal{F}_{L}^{1}(\mathbb{R}^d)$, it is easy to find the two basic facts; the one about the objective value is
\[
f(x_{k+1}) = f\left(y_{k} - s \nabla f(y_{k})\right) \leq f(y_{k}) - \left( s - \frac{Ls^2}{2} \right)\| \nabla f(y_k) \|^2 \leq f(y_{k})
\]
and the other about the gradient is
\[
\nabla f\left(x_{k+1} \right) = \nabla f\left(y_{k} - s \nabla f(y_{k})\right) = \left(I - s\nabla^2 f\left(y_{k} - \theta s \nabla f(y_{k})\right)\right) \nabla f(y_{k}),
\]
where $\theta \in (0,1)$ is a real. Hence, for the convergence rates of the iterative sequence $\{x_k\}_{k=0}^{\infty}$, we conclude them as the following corollary. 

\begin{corollary}
\label{coro:nag-x}
Let $f \in \mathcal{F}_{L}^{1}(\mathbb{R}^d)$. For any step size $0 \leq s \leq 1/L$, the iterates $\{x_{k}\}_{k=0}^{\infty}$ generated by \texttt{NAG} obey
\[
 f(x_{k}) - f(x^\star) \leq \frac{r^2\|x_0 - x^\star\|^2}{2sk(k+r)},\quad \text{and}\quad \lim_{k \rightarrow \infty} \left( k^3 \min_{0\leq i \leq k} \big\| \nabla f(x_{i}) \big\|^2 \right) = 0; 
\]
furthermore, if $r > 2$, then the iterative sequence $\{x_{k}\}_{k=0}^{\infty}$ also satisfy
\[
\lim_{k \rightarrow \infty} \left[ k^2 \left( \min_{0\leq i \leq k}  f(x_{i}) - f(x^\star) \right)  \right] = 0.
\]
\end{corollary}

\section{The implicit-velocity proof}
\label{sec: implicit-proof}

In this section, we investigate the iterative behavior from the implicit-velocity scheme of~\texttt{NAG}~\eqref{eqn: nag-sx} based on the implicit-velocity high-resolution differential equation framework. To analyze the discrete~\texttt{NAG}, the implicit-velocity high-resolution differential equation with its convergence rate is first provided. Then, we demonstrate the convergence rate for~\texttt{NAG} itself from the implicit-velocity scheme~\eqref{eqn: nag-sx}.

\subsection{The implicit-velocity high-resolution differential equation}
\label{subsec: im-hode}

Different from~\citep{shi2021understanding}, we here rearrange the implicit-velocity scheme of~\texttt{NAG}~\eqref{eqn: nag-sx} as
\begin{equation}
\label{eqn: nag-x}
    \frac{x_{k+1} - 2x_{k} + x_{k-1}}{s} + \frac{r+1}{(k+r)\sqrt{s}}\cdot \frac{x_{k}-x_{k-1}}{\sqrt{s}} + \nabla f\left(x_k + \frac{k-1}{k+r}(x_{k}-x_{k-1})\right) = 0.
\end{equation}
Simiarly, plugging the high-resolution Taylor expansion into~\eqref{eqn: nag-x} and taking the $O(\sqrt{s})$-approximation, we can obtain the identical high-resolution differential equation with~\citep{shi2021understanding} as
\[
    \Ddot{X}(t)+\frac{r+1}{t}\Dot{X}(t) + \sqrt{s}\nabla^2 f(X(t))\dot{X}(t) + \left[1+\frac{(r+1)\sqrt{s}}{2t}\right]\nabla f\left(X(t)\right) = 0.
\]
However, to proceed with the proof of~\texttt{NAG} from the implicit-velocity scheme~\eqref{eqn: nag-sx}, we here need to consider the corresponding implicit-velocity high-resolution differential equation as
\begin{equation}
\label{eqn: high-ode}
    \Ddot{X}(t)+\frac{r+1}{t}\Dot{X}(t) + \left[1+\frac{(r+1)\sqrt{s}}{2t}\right]\nabla f\left(X(t)+\sqrt{s}\Dot{X}(t)\right) = 0,
\end{equation}
for any $t \geq (r-1)\sqrt{s}/2$, with $X\left((r-1)\sqrt{s}/2\right) = x_0$ and $\dot{X}\left((r-1)\sqrt{s}/2\right) = - \sqrt{s}\nabla f(x_0)$.
Accompanied by the implicit-velocity high-resolution differential equation~\eqref{eqn: high-ode}, the Lyapunov function is constructed as  
\begin{equation}
\label{eqn: high-ode-lypunov}
    \mathcal{E}(t) = \frac{t(t-r\sqrt{s})}{t-(r+1)\sqrt{s}}\left[t+\frac{(r+1)\sqrt{s}}{2}\right](f(X+\sqrt{s}\Dot{X})-f(x^\ast)) + \frac{1}{2}\big\Vert t\Dot{X}+r(X-x^\ast)\big\Vert^2.
\end{equation}
Then, the convergence rate for the solution  $X = X(t)$ to the implicit-velocity high-resolution differential equation~\eqref{eqn: high-ode} is characterized as the following theorem. 
\begin{theorem}
\label{thm:continuous-implcit-velocity}
Let $f\in \mathcal{F}_1^L(\mathbb{R}^d)$, then the solution $X = X(t)$ to the implicit-velocity high-resolution differential equation~\eqref{eqn: high-ode} satisfies
\begin{equation}
\label{eqn:r=3_continuous}
 f(X+\sqrt{s}\Dot{X})-f(x^\ast) \leq \frac{\mathcal{E}(t_0)}{t^2}, \quad \text{and}\quad  \lim_{t\to\infty} \left(t^3\inf_{t_0\leq u \leq t}\|\nabla f(X+\sqrt{s}\Dot{X})\|^2 \right)=0,
\end{equation}
for all $t\geq t_0=\left(r+2\right)\sqrt{s}$; furthermore, when $r > 2$, the following limitation holds as
\begin{equation}
\label{eqn:r>3_continuous}
    \lim_{t\to\infty} \left[t^2 \left(\inf_{t_0\leq u \leq t}f(X + \sqrt{s}\dot{X}) - f(x^\star) \right)\right] =0.
\end{equation}
\end{theorem}

\begin{proof}
Taking the time derivative of the Lyapunov function $\mathcal{E}(t)$ given in~\eqref{eqn: high-ode-lypunov}, we have
\begin{align*}
\frac{\text{d}\mathcal{E}}{\text{d}t} = & \underbrace{\frac{t(t-r\sqrt{s})}{t-(r+1)\sqrt{s}}\left[t+\frac{(r+1)\sqrt{s}}{2}\right]\left\langle\nabla f(X+\sqrt{s}\Dot{X}), \Dot{X}+\sqrt{s}\Ddot{X}\right\rangle}_{\mathbf{I}} \\
& +\left[2t + \frac{(r+3)\sqrt{s}}{2}-\frac{3(r+1)^2s^{3/2}}{2(t-(r+1)\sqrt{s})^2}\right](f(X+\sqrt{s}\Dot{X})-f(x^\ast)) \\
& - \underbrace{\left\langle t\Dot{X}+r(X-x^\ast), \left[t+\frac{(r+1) \sqrt{s}}{2}\right]\nabla f(X+\sqrt{s}\Dot{X})\right\rangle}_{\mathbf{II}}. 
\end{align*}
With the implicit-velocity high-resolution differential equation~\eqref{eqn: high-ode}, we calculate the first part $\mathbf{I}$ as
\begin{align*}
\mathbf{I} = &\mathrel{\phantom{=}} \underbrace{(t-r\sqrt{s})\left[t+\frac{(r+1)\sqrt{s}}{2}\right] \left\langle\nabla f(X+\sqrt{s}\Dot{X}), \Dot{X}\right\rangle}_{\mathbf{I}_1} \\& -  \frac{ \sqrt{s}\left( t - r\sqrt{s} \right)}{t - (r+1)\sqrt{s}}  \left[ t + \frac{(r+1)\sqrt{s}}{2} \right]^2 \big\| \nabla f(X+\sqrt{s}\Dot{X}) \big\|^2;
\end{align*}
and the second part $\mathbf{II}$ is split into
\begin{align*}
\mathbf{II} = &\mathrel{\phantom{=}}\underbrace{ (t-r\sqrt{s})\left[t+\frac{(r+1)\sqrt{s}}{2}\right] \left\langle\nabla f(X+\sqrt{s}\Dot{X}), \Dot{X}\right\rangle}_{\mathbf{II}_1} \\ & + r\left[t+\frac{(r+1)\sqrt{s}}{2}\right] \left\langle  \nabla f(X+\sqrt{s}\Dot{X}), X + \sqrt{s}\dot{X} - x^\star \right\rangle.
\end{align*}
Obviously, we can find $\mathbf{I}_1 = \mathbf{II}_1$. According to the definition of the convex objective function, we obtain the basic inequality as
\[
    f(X+\sqrt{s}\Dot{X})-f(x^\ast) \leq \left\langle  \nabla f(X+\sqrt{s}\Dot{X}), X + \sqrt{s}\dot{X} - x^\star \right\rangle;
\]
and then estimate the time derivative of the Lyapunov function  $\mathcal{E}(t)$ given in~\eqref{eqn: high-ode-lypunov} as 
\begin{align*}
\frac{\text{d}\mathcal{E}}{\text{d}t} &\leq - \frac{ \sqrt{s}\left( t - r\sqrt{s} \right)}{t - (r+1)\sqrt{s}}  \left[ t + \frac{(r+1)\sqrt{s}}{2} \right]^2 \big\| \nabla f(X+\sqrt{s}\Dot{X}) \big\|^2 \\
                                      & \mathrel{\phantom{\leq}} - \left[ (r-2)t + \frac{(r^2-3)\sqrt{s}}{2} \right]\left(f(X+\sqrt{s}\Dot{X})-f(x^\ast)\right) \\
                                      & \leq -  \sqrt{s} t^2 \big\| f(X+\sqrt{s}\Dot{X}) \big\|^2 - (r-2)t \left(f(X+\sqrt{s}\Dot{X})-f(x^\ast)\right).
\end{align*}
Therefore, the Lyapunov function $\mathcal{E}(t)$ given in~\eqref{eqn: high-ode-lypunov} is decreasing with time $t \geq t_0$. Meanwhile, the following integral inequality can be derived as
\[
 \int_{t_0}^{t}\left[ \sqrt{s}u^2 \big\| \nabla f(X+\sqrt{s}\Dot{X}) \big\|^2 + (r-2)u \left(f(X+\sqrt{s}\Dot{X})-f(x^\ast)\right)\right] du \leq \mathcal{E}(t_0).
\] 
The limitation~\eqref{eqn:r=3_continuous} hence is derivation as
\[
0 \leq \frac{7\sqrt{s}}{24} \lim_{t \rightarrow \infty} \left( t^3 \inf_{t_0\leq u \leq t} \big\| \nabla f(X+\sqrt{s}\Dot{X})) \big\|^2 \right) \leq \lim_{t\rightarrow \infty}\int_{\frac{t}{2}}^{t}u^2\big\| \nabla f(X+\sqrt{s}\Dot{X}) \big\|^2 du= 0.
\]
Furthermore, when $r > 2$, we complete the derivation~\eqref{eqn:r>3_continuous} as
\[
0  \leq \frac{3}{8} \lim_{k \rightarrow \infty} \left[ t^2 \left( \inf_{t_0\leq u \leq t}  f(X+\sqrt{s}\Dot{X}) - f(x^\star) \right)  \right] 
  \leq  \int_{\frac{t}{2}}^{t} u \left(f(X+\sqrt{s}\Dot{X})-f(x^\ast)\right) du = 0.
\]
\end{proof}

\subsection{The implicit-velocity scheme}
\label{subsec: implicit-velocity-scheme}

To consider the phase-space representation of the implicit-velocity scheme of~\texttt{NAG}~\eqref{eqn: nag-sx}, different from the phase-space representation of the gradient-correction scheme of~\texttt{NAG}~\eqref{eqn: nag-sy}, we here takes the velocity iterates as $v_{k} = (x_k - x_{k-1})/\sqrt{s}$, which is the implicit scheme. Then the phase-space representation of the implicit-velocity scheme of~\texttt{NAG}~\eqref{eqn: nag-sx} can be written as
\begin{equation}
\label{eqn: phase-x}
\left\{\begin{aligned}
& x_{k+1} - x_{k} = \sqrt{s}v_{k+1}, \\
&v_{k+1} - v_{k} = -\frac{r+1}{k+r}v_k - \sqrt{s}\nabla f\left(x_k + \frac{k-1}{k+r}\sqrt{s}v_k\right).
\end{aligned} \right.
\end{equation}
Thus, we call~\eqref{eqn: phase-x} as~\textbf{implicit-velocity phase-space representation of~\texttt{NAG}}. Different from the phase-space  representation of gradient-correction scheme~\eqref{eqn:phase-nag}, here the first line for the position sequence $\{x_k\}_{k=0}^{\infty}$ is an implicit scheme, while the second line for the velocity sequence $\{v_{k}\}_{k=0}^{\infty}$ is an explicit scheme.  Recall the~\texttt{NAG}; the second scheme can be written as
\begin{equation}
\label{eqn: nag-2-yx}
    y_k = x_k + \frac{k-1}{k+r}\cdot\sqrt{s}v_k.
\end{equation}
Therefore, the implicit-explicit implicit-velocity phase-space representation of~\texttt{NAG}~\eqref{eqn: phase-x} can be expressed as 
\begin{equation}
\label{eqn: phase-xy}
\left\{\begin{aligned}
& x_{k+1} - x_{k} = \sqrt{s}v_{k+1}, \\
&v_{k+1} - v_{k} = -\frac{r+1}{k+r}v_k - \sqrt{s}\nabla f\left(y_k\right).
\end{aligned} \right.
\end{equation}
Corresponding to the implicit-explicit implicit-velocity phase-space representation of~\texttt{NAG}~\eqref{eqn: phase-x}, the convergence rates are calculated by the new  Lyapunov function constructed as
\begin{equation}\label{eqn: lyapunov-implicit}
    \mathcal{E}(k) = sk(k+r)\left(f\left(y_{k-1}\right)-f(x^\ast)\right) + \frac{1}{2}\left\|\sqrt{s}(k-1)v_k+r(x_k-x^\ast)\right\|^2,
\end{equation}
The relationship with the previous Lyapunov function~\eqref{eqn: lypunov-positive} will be discussed in~\Cref{sec: discussion}.


\begin{proof}[Proof of~\Cref{thm:nag} in implicit-velocity scheme]
The second equality of the implicit-explicit implicit-velocity phase-space representation~\eqref{eqn: phase-xy} can be equivalently written as
\begin{equation}
\label{eqn: phase-xy-1}
    (k+r)v_{k+1} - (k-1)v_k = -(k+r)\sqrt{s}\nabla f\left(y_k\right).
\end{equation}
Here, we also separate it into three steps to clarify the proof.
\begin{itemize}
\item[\textbf{(1)}] The difference between the second terms of Lyapunov function~\eqref{eqn: lyapunov-implicit}, $\sqrt{s}kv_{k+1}+r(x_{k+1}-x^\ast)$ and $\sqrt{s}(k-1)v_k+r(x_k-x^\ast)$, is 
\begin{multline*}
 \left[\sqrt{s}kv_{k+1}+r(x_{k+1}-x^\ast)\right]  - \left[\sqrt{s}(k-1)v_k+r(x_k-x^\ast)\right] \\
    =  \sqrt{s}kv_{k+1} - \sqrt{s}(k-1)v_{k} + r(x_{k+1}-x_{k}). 
\end{multline*}
With the implicit-explicit implicit-velocity phase-space representation~\eqref{eqn: phase-xy} and its equivalent form~\eqref{eqn: phase-xy-1}, the difference can be computed as
\begin{align*}
 \left[\sqrt{s}kv_{k+1}+r(x_{k+1}-x^\ast)\right]  - \left[\sqrt{s}(k-1)v_k+r(x_k-x^\ast)\right]   & =  \sqrt{s}((k+r)v_{k+1}-(k-1)v_{k}) \\
    &= -s(k+r)\nabla f(y_{k}).
\end{align*}
Similarly, with~\eqref{eqn: nag-2-yx} and~\eqref{eqn: phase-xy}, we can calculate the difference between $y_k$ and $y_{k-1}$ as 
\begin{align*}
y_k - y_{k-1} & = \left[x_k + \frac{k-1}{k+r} \cdot \sqrt{s}v_k\right] - \left[x_{k-1} + \frac{k-2}{k+r-1}\cdot \sqrt{s}v_{k-1}\right] \\
              & = \frac{k-1}{k+r}\cdot \sqrt{s}v_k + \sqrt{s}\left( v_{k} - \frac{k-2}{k+r-1} \cdot v_{k-1} \right) \\
              & = \frac{k-1}{k+r}\cdot \sqrt{s}v_k - s \nabla f(y_{k-1}).     
\end{align*}

\item[\textbf{(2)}]
Then, with the subscripts labeled, we calculate the difference between $\mathcal{E}(k+1)$ and $\mathcal{E}(k)$ as
\begin{align*}
\mathcal{E}(k+1) -\mathcal{E}(k) &=  sk(k+r)\left(f\left(y_{k}\right) - f\left(y_{k-1}\right)\right) + s(2k+r+1)\left(f\left(y_{k}\right) - f(x^\ast)\right) \\
& \mathrel{\phantom{=}} -  \left\langle s(k+r)\nabla f(y_k), \sqrt{s}(k-1)v_k + r(x_k - x^\star)\right\rangle   + \frac{s^2(k+r)^2}{2} \| \nabla f(y_k) \|^2 \\ 
 &=  sk(k+r)\left(f\left(y_{k}\right) - f\left(y_{k-1}\right)\right) + s(2k+r+1)\left(f\left(y_{k}\right) - f(x^\ast)\right) \\
& \mathrel{\phantom{=}} -  \left\langle s(k+r)\nabla f(y_k), \frac{\sqrt{s}k(k-1)v_k}{k+r} + r(y_k - x^\star)\right\rangle   + \frac{s^2(k+r)^2}{2} \| \nabla f(y_k) \|^2 \\ 
&=  \underbrace{sk(k+r)\left(f\left(y_{k}\right) - f\left(y_{k-1}\right)\right)}_{\mathbf{I}_1} + \underbrace{s(2k+r+1)\left(f\left(y_{k}\right) - f(x^\ast)\right)}_{\mathbf{I}_2} \\
& \mathrel{\phantom{=}} -\underbrace{ s^{\frac{3}{2}}k(k-1) \left\langle \nabla f(y_k), v_k \right\rangle}_{\mathbf{II}_1} - \underbrace{sr(k+r) \left\langle \nabla f(y_k), y_k - x^\star\right\rangle}_{\mathbf{II}_2} \\ & \mathrel{\phantom{=}}  + \frac{s^2(k+r)^2}{2} \| \nabla f(y_k) \|^2,
\end{align*}
where the second equality follows from the second equality of~\texttt{NAG}~\eqref{eqn: nag-2-yx}.

\item[\textbf{(3)}] For any $f \in \mathcal{F}_{L}^{1}$, we can estimate $\mathbf{I}_1$ as 
\begin{align*}
 \mathbf{I}_1 & = sk(k+r) \left( f\left(y_{k}\right) - f\left(y_{k-1}\right) \right) \\
              & \leq \left\langle\nabla f\left(y_{k}\right), y_k-y_{k-1}\right\rangle -\frac{sk(k+r)}{2L}\left\|\nabla f\left(y_{k}\right) - \nabla f\left(y_{k-1}\right)\right\|^2 \\
             & = s^{\frac{3}{2}}k(k-1)\left\langle\nabla f\left(y_{k}\right), v_{k}\right\rangle + s^2k(k+r)\left\langle\nabla f\left(y_{k}\right), \nabla f\left(y_{k-1}\right)\right\rangle -\frac{sk(k+r)}{2L}\left\|\nabla f\left(y_{k}\right) - \nabla f\left(y_{k-1}\right)\right\|^2,
\end{align*}
which leads to the difference between $\mathbf{I}_1$ and $\mathbf{II}_1$ as
\[
\mathbf{I}_1 -  \mathbf{II}_1 = s^2k(k+r)\left\langle\nabla f\left(y_{k}\right), \nabla f\left(y_{k-1}\right)\right\rangle -\frac{sk(k+r)}{2L}\left\|\nabla f\left(y_{k}\right) - \nabla f\left(y_{k-1}\right)\right\|^2.
\]

Similarly, the difference between $\mathbf{I}_2$ and $\mathbf{II}_2$ can be computed as 
\begin{align*}
\mathbf{I}_2 - \mathbf{II}_2 & = s (2k+r+1) \left(   f\left(y_{k}\right) - f(x^\ast) \right) - sr(k+r) \left\langle \nabla f(y_k), y_k - x^\star\right\rangle \\
             & \leq -s \left[(r-2)k+r^2 -r-1\right]\left(   f\left(y_{k}\right) - f(x^\ast) \right) -\frac{s r(k+r)}{2L}\left\|\nabla f\left(y_{k}\right)\right\|^2.
\end{align*}
Hence, we can obtain the iterative difference of the Lyapunov function~\eqref{eqn: lyapunov-implicit} as
\begin{align*}
\mathcal{E}(k+1) - \mathcal{E}(k) & = - s^{2}k(k+r)\langle \nabla f(y_{k}), \nabla f(y_{k-1})\rangle - \frac{sk(k+r)}{2L} \| \nabla f(y_{k}) - \nabla f(y_{k-1}) \|^2  \\
                                  & \mathrel{\phantom{=}} -\left[(r-2)k+r^2 -r-1\right]\left(   f\left(y_{k}\right) - f(x^\ast) \right) -\frac{s r(k+r)}{2L}\left\|\nabla f\left(y_{k}\right)\right\|^2\\
                                  & \mathrel{\phantom{=}} + \frac{s^2(k+r)^2}{2} \| \nabla f(y_k) \|^2.
\end{align*}
Furthermore, when $0 < s \leq 1/L $, the difference between $\mathcal{E}(k+1)$ and $\mathcal{E}(k)$ satisfies 
\[
\mathcal{E}(k+1) - \mathcal{E}(k) \leq - \frac{s^2k(k+r)}{2}\big\| \nabla f(y_{k-1}) \big\|^2 - s \left[ (r-2)k + r^2 - r - 1  \right]\left( f(y_{k}) - f(x^\star) \right),
\]
which is exactly the same as~\eqref{eqn:lyapunov-monotone}. Therefore, the lines following~\eqref{eqn:lyapunov-monotone} can be used to complete the proof of~\Cref{thm:nag}.
\end{itemize}
\end{proof}

\section{Discussion}
\label{sec: discussion}

In this study, we refomulate and simplify the Lyapunov function for~\texttt{NAG} constructed in~\citep{shi2021understanding}, which is given in~\eqref{eqn: lypunov-positive} as
\[
\mathcal{E}(k) = sk(k+r)\left( f(y_{k-1}) - f(x^\star) \right) + \frac12 \| \sqrt{s}kv_{k-1} + r(y_k - x^\star) + sk\nabla f(y_{k-1}) \|^2;
\]
meanwhile, we propose a new Lyapunov function for the implicit-velocity scheme in~\eqref{eqn: lyapunov-implicit} as
\[
  \mathcal{E}(k) = sk(k+r)\left(f\left(y_{k-1}\right)-f(x^\ast)\right) + \frac{1}{2}\left\|\sqrt{s}(k-1)v_k+r(x_k-x^\ast)\right\|^2,
\]
which does not include the gradient term. Through their corresponding phase-space representations,~\eqref{eqn:phase-nag} and~\eqref{eqn: phase-x}, both the Lyapunov functions~\eqref{eqn: lypunov-positive} and~\eqref{eqn: lyapunov-implicit} are identical and equal to
\begin{equation}
\label{eqn: lyapunov-new}
\mathcal{E}(k) = sk(k+r)\left( f(y_{k-1}) - f(x^\star) \right) + \frac12 \| k (y_k - x_k) + r(y_k - x^\star) \|^2,
\end{equation}
where the essential difference from the Lyapunov function constructed in~\citep[Theorem 6]{su2016differential} is that the potential function here is $ f(y_{k-1}) - f(x^\star)$ other than  $ f(x_{k}) - f(x^\star)$, although the coefficients of the potential function are also different. Recall in~\citep[Theorem 6]{su2016differential} that the difference in the potential function  is
\begin{align*}
& s(k+r)^2\left( f(x_{k+1}) - f(x^\star) \right) - s(k+r-1)^2\left( f(x_{k}) - f(x^\star) \right) \\
\leq &s(k+r) \left[ (k+r)\left( f(x_{k+1}) - f(x^\star) \right) - k \left( f(x_{k}) - f(x^\star) \right)\right] - \left[ (r-2)k + (r-1)^2\right]\left( f(x_{k}) - f(x^\star) \right) \\
= &s(k+r) \left[ k\left( f(x_{k+1}) - f(x_k) \right) + r \left( f(x_{k+1}) - f(x^\star) \right)\right] - \left[ (r-2)k + (r-1)^2\right]\left( f(x_{k}) - f(x^\star) \right) \\
= &\underbrace{sk(k+r)\left( f(x_{k+1}) - f(x_k) \right)}_{\mathbf{J}_1} + \underbrace{sr(k+r) \left( f(x_{k+1}) - f(x^\star) \right)}_{\mathbf{J}_2} - \left[ (r-2)k + (r-1)^2\right]\left( f(x_{k}) - f(x^\star) \right). 
\end{align*}
The difference in the potential function shown in~\Cref{subsec: implicit-velocity-scheme} is
\begin{multline*}
s(k+1)(k+r+1)\left( f(y_{k}) - f(x^\star) \right) - sk(k+r)\left( f(y_{k-1}) - f(x^\star) \right) \\
\leq \underbrace{sk(k+r)\left( f(y_{k}) - f(y_{k-1}) \right)}_{\mathbf{I}_1} + \underbrace{s(2k+r+1)\left( f(y_{k}) - f(x^\star) \right)}_{\mathbf{I}_2}. 
\end{multline*}
Indeed, we here need to use $\mathbf{J}_1$ and $\mathbf{J}_2$ instead of $\mathbf{I}_1$ and $\mathbf{I}_2$, respectively, to calculate the difference between the Lyapunov functions, $\mathcal{E}(k+1)$ and $\mathcal{E}(k)$. Taking the inequality for the convex function in~\citep[(22)]{su2016differential}, $\mathbf{J}_1$ and $\mathbf{J}_2$ are estimated as
\begin{align*}
  f(x_{k+1}) - f(x_k) & \leq \left\langle \nabla f(y_k), y_{k} - x_{k}  \right\rangle - \frac{s}{2} \| \nabla f(y_k) \|^2 \\
                      & =    \left\langle \nabla f(y_k), y_{k} - y_{k-1}  \right\rangle + s\left\langle \nabla f(y_k), \nabla f(y_{k-1})  \right\rangle   - \frac{s}{2} \| \nabla f(y_k) \|^2
\end{align*}
and
\begin{align*}
  f(x_{k+1}) - f(x^\star)  \leq \left\langle \nabla f(y_k), y_{k} - x^\star  \right\rangle - \frac{s}{2} \| \nabla f(y_k) \|^2, \mathrel{\phantom{   - \frac{s}{2} \| \nabla f(y_k) \|^2fdsfsdfsd}}
\end{align*}
where we find the cross term $s\left\langle \nabla f(y_k), \nabla f(y_{k-1})  \right\rangle$ appears in the estimate of $\mathbf{J}_1$. In other words, the convex inequality in~\citep[(22)]{su2016differential} is not tight, which directly leads to 
\[
\mathcal{E}(k+1) - \mathcal{E}(k) \leq - \left[ (r-2)k + (r-1)^2\right]\left( f(x_{k}) - f(x^\star) \right) -\frac{s r(k+r)}{2L}\left\|\nabla f\left(y_{k}\right)\right\|^2.
\]
Hence, we cannot find the gradient norm acceleration since the inequality~\citep[(22)]{su2016differential} is not tight.

For a Hamilton system, there exist two kinds of symplectic (semi-implicit or semi-explicit) Euler schemes, an explicit-implicit scheme and an implicit-explicit scheme, which share the same property~\citep{haier2006geometric}. Similarly, we can say here for the Lyapunov function~\eqref{eqn: lyapunov-new}, there exist two kinds of phase-space representations of~\texttt{NAG}, the explicit-implicit gradient-correction scheme~\eqref{eqn:phase-nag} and the implicit-explicit implicit-velocity scheme~\eqref{eqn: phase-x}, respectively. Different from the case of the Hamilton system, there exists a striking distinctness of two kinds of phase-space representations for the nonlinear convergence of~\texttt{NAG}. Still, we find the two kinds of phase-space representations are substantially equivalent. Furthermore, for the explicit-implicit gradient-correction scheme~\eqref{eqn:phase-nag}, we use the Lyapunov function with the form~\eqref{eqn: lypunov-positive} to characterize the convergence rate, where the gradient term in the mixed term of~\eqref{eqn: lypunov-positive} corresponds to the gradient correction; while for the implicit-explicit gradient-correction scheme~\eqref{eqn: phase-x}, the Lyapunov function can be rewritten as
\begin{equation}
\label{eqn: lyapunov-new1}
\mathcal{E}(k) = sk(k + r)\left[f\left( x_{k-1} + \frac{(k-2)\sqrt{s} v_{k - 1}}{k+r-1}  \right) - f(x^\star) \right]  + \frac{1}{2}\left\|\sqrt{s}(k-1)v_k+r(x_k-x^\ast)\right\|^2,
\end{equation}
where the velocity is implicitly included in the potential and the gradient term disappears in the mixed term from the aspect of the iterative sequence $\{ x_k \}_{k=0}^{\infty}$. Thus, we answer the question proposed in~\Cref{sec: intro}, the implicit-explicit implicit-velocity scheme of~\texttt{NAG}~\eqref{eqn: nag-sx} shares the same convergence rate with the explicit-implicit gradient-correction scheme~\eqref{eqn: nag-sy}.

Finally, we remark on some possible extensions. In this study, we simplify the proof from the explicit-implicit gradient-correction scheme in~\citep{shi2021understanding} and discuss the convergence behavior of~\texttt{NAG} from the implicit-explicit implicit-velocity scheme, or the iterative sequence $\{x_k\}_{k=0}^{\infty}$. Meanwhile, we find the implicit-explicit implicit-velocity scheme is equivalent to the explicit-implicit gradient-correction scheme based on the identical Lyapunov function to derive the convergence rate. Still, extending the implicit-velocity high-resolution differential equation framework looks more accessible beyond smooth convex optimization in the Euclidean setting. Moreover, an accelerated algorithm for mirror descent proposed in~\citep{krichene2015accelerated}, which is on the grounds of low-resolution differential equations, is not available in practice. Hence, it is exciting to provide an accelerated mirror descent method available in practice using the high-resolution differential equation framework in non-Euclidean spaces. Recently, some appealing research in statistics relates to gradient-based algorithms, such as using gradient descent to investigate the implicit sparsity-inducing mechanism~\citep{jordan2021self}. It is also exciting to explore the implicit sparsity-inducing mechanism based on the high-resolution differential equation framework.

{\small
\bibliographystyle{abbrvnat}
\bibliography{reference}
}
\end{document}